\documentclass[sn-mathphys,Numbered]{sn-jnl}% Math and Physical Sciences Reference Style
\usepackage[arrow,matrix]{xy}
\usepackage{graphicx}%
\usepackage{multirow}%
\usepackage{amsmath,amssymb,amsfonts}%
\usepackage{amsthm}%
\usepackage{mathrsfs}%
\usepackage[title]{appendix}%
\usepackage{xcolor}%
\usepackage{textcomp}%
\usepackage{manyfoot}%
\usepackage{booktabs}%
\usepackage{algorithm}%
\usepackage{algorithmicx}%
\usepackage{algpseudocode}%
\usepackage{listings}%
\theoremstyle{thmstyleone}%
\theoremstyle{thmstyletwo}%

\theoremstyle{thmstylethree}%
\newtheorem{thm}{Theorem}

\newtheorem{defn}{Definition}
\newtheorem{lem}{Lemma}
\newcommand{\Res}{\operatornamewithlimits{Res}}
\newcommand{\dslash}{/\!\!/}

\newcommand{\mk}{\mathfrak}
\newcommand{\mc}{\mathcal}
\newcommand{\mb}{\mathbb}

\newcommand{\GL}{\operatorname{GL}}

\newcommand{\mo}{\mathcal{O}}
\newcommand{\E}{\mathcal{E}}

\newcommand{\PP}{\mathbb{P}}

\def\Eu{\mathop{\mathrm{Eu}}\nolimits}

\def\exc{\mathop{\mathrm{exc}}\nolimits}
\def\pt{\mathop{\mathrm{pt}}\nolimits}

\def\C{\mathbb{C}}

\def\Z{\mathbb{Z}}

\def\Hom{\mathop{\mathrm{Hom}}\nolimits}

\def\rk{\mathop{\mathrm{rk}}\nolimits}

\def\pt{\mathop{\mathrm{pt}}\nolimits}
\def\kloc{K^{\text{loc}}}

\begin{document}

\title[Residue formula for flag manifold of type $A$ from wall-crossing]
{Residue formula for flag manifold of type $A$ from wall-crossing}

\author*[1,2]{\fnm{Ryo} \sur{Ohkawa}}
\email{ohkawa.ryo@gmail.com}
 
 \affil*[1]{\orgdiv{Osaka Central Advanced Mathematical Institute}, 
\orgname{Osaka Metropolitan University}, 
\orgaddress{
%\street{Street}, \city{Osaka-si}, 
\postcode{558-8585}, \state{Osaka}, \country{Japan}}}

\affil*[2]{\orgdiv{Research Institute for Mathematical Sciences}, \orgname{Kyoto University}, 
\orgaddress{
%\street{Oiwake-cho}, 
%\city{Kyoto-si}, 
\postcode{606-8502}, \state{Kyoto}, \country{Japan}}}

\abstract{
We consider equivariant integrals on flag manifolds of 
type $A$.
%, especially Grassmannian manifolds. 
Using a computational method inspired by the theory of 
wall-crossing formulas by Takuro Mochizuki, we re-prove 
residue formulas for equivariant integrals given by 
Weber and Zielenkiewicz. 
As an application, we give the determinantal formula of 
the Grothendieck polynomial by properly setting 
$K$ theory classes.
}


%%\pacs[JEL Classification]{D8, H51}

%%\pacs[MSC Classification]{35A01, 65L10, 65L12, 65L20, 65L70}

\maketitle

\section{Introduction}
\label{sec1}

We study equivariant integrals over flag manifolds of type $A$ via wall-crossing phenomena
and give explicit formula.
This gives another proof of residue formula proved by Weber-Zielenkiewicz \cite{WZ}
using torus equivariant integrals over $K$-theory classes on Grassmann manifolds.

In \cite{WZ}, they study a homogeneous space
with the maximal torus action.
One of the examples is the Grassmannian $Gr(r, d)$ of the $d$-dimensional spaces in $\C^{r}$.
The Grassmannian $Gr(r,d)$ is a GIT quotient $\Hom(\C^{d}, \C^{r}) \dslash \GL(\C^{d})$, 
that is the geometric quotient of the open set consisting of the injective maps $\C^{d} \to \C^{r}$.
A generalization of the theory of Jeffrey-Kirwan \cite{JK} or Guillemin-Kalkman \cite{GK} applied 
to the homogeneous spaces leads to explicit residue formulas for the push-forward in equivariant cohomology. 

The Jeffrey-Kirwan method can be applied for a class 
of generalized cohomology theories as shown in \cite{Me}. 
For the Grassmannian $Gr(r, d)$ the formula for the push-forward obtained via symplectic reduction is the following
\[
\Pi_{\ast}(\E)
={1 \over d!}
\Res_{u_{1}, \ldots, u_{d}=0, \infty} 
{
f(u_{1}, \ldots, u_{d}) \prod_{i \neq j} 
(1-u_{i}/u_{j})
\over 
\prod_{\alpha=1}^{r} \prod_{i=1}^{d} 
(1- u_{i}/x_{\alpha})}
{du_{1} \cdots du_{d}
\over 
u_{1} \cdots u_{d}}
\]
where $\E$ is vector bundle defined from $\GL(\C^{d})$-representation $E$, and 
$f$ is the character when $E$ is restricted to a maximal torus. 
The operation $\Res_{u_{1}, \ldots, u_{d}=0,\infty}$ is the iterated residue 
$\displaystyle  \Res_{u_{1}=0, \infty} \circ \Res_{u_{2}=0, \infty} \circ \cdots
\circ \Res_{u_{d}=0, \infty}$ where 
the symbol $\displaystyle \Res_{u = 0,\infty}$ denotes 
taking a sum of residues at $0$ and $\infty$.
In our case the order of taking the residues does not 
matter.

In this paper, we study integrals over flag manifolds of type $A$ 
using computational method inspired by the theory of wall-crossing formulas developed by Takuro Mochizuki \cite{Mo}.
For example, we consider variation of stability parameters for $\GL(\C^{d})$-action on $\Hom(\C^{d}, \C^{r})$
corresponding to characters $\det^{\pm 1} \colon \GL(\C^{d}) \to \C^{\ast}, g \mapsto \det(g)^{\pm 1}$.
When we take an ample $\GL(\C^{d})$-equivariant line bundle $\mo_{\Hom(\C^{d}, \C^{r})} \otimes \C_{\det^{1}}$,
we have $Gr(r,d)=\Hom(\C^{d}, \C^{r}) \dslash_{\det^{1}} \GL(\C^{d})$ as a GIT quotient.
On the other hand, we have $\emptyset=\Hom(\C^{d}, \C^{r}) \dslash_{\det^{-1}} \GL(\C^{d})$. 
This leads to wall-crossing phenomena between $Gr(r, d)$ and a empty set.
Since wall-crossing terms are given by integrals over $Gr(r, d')$ for $d' < d$, we get a recursion formula
for integrals over $Gr(r, d)$.

To this end, we introduce a {\it master space} $\mc M$ 
explained below in a slightly generalized situation.
Method using the space $\mc M$ is motivated 
by the theory developed by Mochizuki \cite{Mo}.
We consider an algebraic torus $\mb T=(\C^{\ast})^{s}$
and $\mb T$-equivariant vector bundle $\mc E$ over a 
base $\mb T$-manifold $B$.
Then $\mc M$ is a projective bundle over a flag bundle
$Fl_{B}(\mc E;d-1)$ over $B$ (cf. \S \ref{subsec:flag}).
Localization technique applied to $\mc M$ allows us to 
reduce $K$-theoretic integrals over $Fl_{B}(\mc E; d)$
to ones over $Fl_{B}(\mc E; d-1)$.
As a main result, we have Theorem \ref{thm:main}. 

As an application of the main result, we give 
determinantal formula Theorem \ref{thm:groth} for the 
{\it refined shifted factorial Grothendieck polynomial}
defined in Definition \ref{rsfg}. 
Determinantal formulas for shifted factorial Grothendieck polynomial using Schubert calculus are presented in \cite{HIMN}
generalizing the Jacobi-Trudi formula for Schur polynomial.
Another generalization is obtained in \cite{HJKSS}.
We further generalize this formula using the residue formula.
These kind of determinantal formulas are also studied using various methods in \cite{I}, \cite{MS}.

In \cite{WZ} and \cite{Z}, the residue formulas for 
more general types of flag manifoles are obtained.
It is natural to ask whether it is possible or not to 
apply our method to other types of flag manifolds.
% as studied in \cite{WZ}.
See \cite{M} and \cite{MOSWW} for algebro-geometric 
descriptions of general types of flag manifolds.
It is also interesting to try to apply our method to 
push-forwards in quantum $K$-theory as studied 
in \cite{KLNS}.
Push-forward formulas for the Plucker classes 
in cohomology are also obtained in 
\cite{KT}.

This paper is organized as follows.
In Section 2, we recall $\mb T$-equivariant $K$-theory, 
and introduce master spaces $\mc M$ 
for flag bundles of type $A$ and a
torus $\C^{\ast}_{u}$ acting on $\mc M$.
The localization theorem for this torus 
$\C^{\ast}_{u}$-action gives us 
equations between push-forwards from 
$Fl_{B}(\mc E;d)$ and $Fl_{B}(\mc E;d-1)$. 
In section 3, we apply this equations repeatedly
to get the wall-crossing formula Theorem \ref{thm:main}
for $K$-theoretic push-forwards.
We also consider cohomological push-forwards and 
deduce Theorem \ref{thm:cohom}.
In section 4, we apply these formula to compute 
integrals of explicit $K$-theoretic classes to get
the determinantal formula Theorem \ref{thm:groth} 
for the 
refined shifted factorial Grothendieck polynomial.

The author is grateful for Hidetoshi Awata, Takeshi 
Ikeda, 
Koji Hasegawa, Ayumu Hoshino, Hiroaki Kanno, Hitoshi 
Konno, Takuro Mochizuki, Kohei Motegi, Satoshi Naito, 
Hiraku Nakajima, Hiroshi Naruse, Masatoshi Noumi, 
Takuya Okuda,  
Yusuke Ohkubo, Yoshihisa Saito and Jun'ichi Shiraishi 
for 
discussion.
He has had the generous support and encouragement of 
Masa-Hiko Saito.
The author would like to appreciate referees for careful 
suggestions and corrections.
They made the manuscript greatly improved.
In particular, our method became applied to more general 
setting for push-forwards from flag bundles thanks
to the referee's suggestion. 

%%%%%%%%%%%%%%%%%%%%%%%%%%%%%%%%%%%%%%%%%%%%%%%%%
%%%%%%%%%%%%%%%%%%%%%%%%%%%%%%%%%%%%%%%%%%%%%%%%%
%%%%%%%%%%%%%%%%%%%%%%%%%%%%%%%%%%%%%%%%%%%%%%%%%
\section{Master space}
\label{sec:master}

\subsection{Localization}

We consider an algebraic torus $\mb T=(\C^{\ast})^{s}$.
For a character $\chi \colon \mb T \to \C^{\ast}$, 
the weight space $\C_{\chi}$ denotes the one-dimensional
$\mb T$-representation $\C$ with the $\mb T$-action
defined by $t \cdot 1 = \chi(t)$ for $t \in \mb T$.
For a coordinate $(x_{1}, \ldots, x_{s})$ of $\mb T$, 
we write by $x_{1}^{m_{1}} \cdots x_{s}^{m_{s}}$ the 
character sending $t \in \mb T$ to 
$x_{1}(t)^{m_{1}} \cdots x_{s}(t)^{m_{s}} \in \C^{\ast}$.
The symbol $x_{1}^{m_{1}} \cdots x_{s}^{m_{s}}$ also
denotes the weight space $\C_{x_{1}^{m_{1}} \cdots
x_{s}^{m_{s}}}$.

For a smooth algebraic variety $M$ with the 
$\mb T$-action, 
we consider the $K$-theory $K_{\mb T}(M)$ of 
$\mb T$-equivariant vector bundles on $M$.
We can regard $K_{\mb T}(M)$ as the Grothendieck ring 
of the category of $\mb T$-equivariant 
{\it locally free sheaves} on $M$.
The multiplication is induced by the tensor product.
For a morphism $M_{1} \to M_{2}$ of algebraic varieties
and a vector bundle $\mc E$ on $M_{2}$, 
we often write the pull-back $\mc E|_{M_{1}}$ 
by the same letter $\mc E$.
Since we always assume that $M$ is smooth in
this paper, we can identify 
$K_{\mb T}(M)$ with the Grothendieck group of the
category of $\mb T$-equivariant 
{\it coherent sheaves} on $M$ by 
\cite[Proposition 2.1]{N}. 
In particular, when $M$ is equal to the one point pt,
we have the identification $K_{\mb T}(\pt)=
\Z[x_{1}^{\pm 1}, \ldots, x_{s}^{\pm 1}]$ with the 
representation ring by
the weight space decomposition.

For a $\mb T$-equivariant vector bundle 
$\mc E$ over $M$, 
we write by $\mc E^{\vee}$ and $\wedge^{k} \mc E$ 
the dual
of $\mc E$ and the $k$-th exterior power 
of $\mc E$ respectively.
For two $\mb T$-equivariant
vector bundles $\mc E$ and $\mc F$ over $M$, we 
set $\mc Hom(\mc E, \mc F)=\mc E^{\vee} \otimes \mc F$.
When $\rk \mc E=r$, we set 
%$\det \mc E = \wedge^{r} \mc E$, and 
$\wedge_{-1} \mc E= \sum_{i=0}^{r} 
(-1)^{i} \cdot \wedge^{i} \mc E$.
We consider a multiplicative subset
\[
S_{\mb T, M}=
\lbrace 
\wedge_{-1} \mc E \mid  
\mc E \colon \mb T \text{-equivariant vector bundle
such that } \wedge_{-1} \mc E \neq 0
\rbrace \subset K_{\mb T}(M),
\]
and the localized ring 
$\kloc_{\mb T} (M)
= S_{\mb T, M}^{-1} K_{\mb T}(M)$ with respect
to $S_{\mb T, M}$.

For a polynomial 
$g(Y_{1}, \ldots, Y_{r}) \in \kloc_{\mb T}(M)
[Y_{1}, \ldots, Y_{r}]$ and 
$\mb T$-equivariant vector bundle $\mc E$, we set 
\[
g(\wedge^{\bullet} \mc E) = 
g(\wedge^{1} \mc E, \ldots, \wedge^{r} \mc E)
\in \kloc_{\mb T}(M).
\]
In particular, when $g(Y_{1}, \ldots, Y_{r}) = 
\sum_{i=1}^{r}(-1)^{i}Y_{i}$, 
we have
$g(\wedge^{\bullet} \mc E)=\wedge_{-1} \mc E$.
For a symmetric polynomial 
$f(x_{1}, \ldots, x_{s})$, we can write 
\[
f(x_{1}, \ldots, x_{s}) = 
g(e_{1}(x_{1}, \ldots, x_{s}), \ldots, 
e_{s}(x_{1}, \ldots, x_{s}) )
\]
using a polynomial $g(Y_{1}, \ldots, Y_{s})$ 
and the elementary symmetric polynomials 
$e_{1}(x_{1}, \ldots, x_{s}), \ldots, 
e_{s}(x_{1}, \ldots, x_{s})$.
Then we set $f(\mc E) = 
g(\wedge^{\bullet}\E)$.

For $\mb T$-equivariant vector bundles $\mc E_{1}, 
\mc E_{2}$ on $M$, the {\it virtual vector bundle}
$\mc E_{1} - \mc E_{2}$ can be regarded as an element in 
$K_{\mb T}(M)$, and all elements of $K_{\mb T}(M)$
are presented in this way.
If it is necessary, we often replace $\mb T$ with 
$\mb T \times \C^{\ast}_{u}$ and consider elements 
\[
%\wedge_{-1} \colon K_{\mb T}(M) \to \kloc_{\mb T} (M),
%\quad  
%\mc E_{1} - \mc E_{2} \mapsto 
%\wedge_{-u} \left( \mc E_{1} - \mc E_{2} \right)
%=
{ \wedge_{-1} ( \mc E_{1} \otimes \C_{u} )
\over 
\wedge_{-1} ( \mc E_{2} \otimes \C_{u} )
}
=
{ \sum_{i} (-u)^{i} \wedge^{i} \mc E_{1}  
\over 
\sum_{i} (-u)^{i} \wedge^{i} \mc E_{2}
} 
\]
in $\kloc_{\mb T \times \C^{\ast}_{u}} (M) = \kloc_{\mb T} (M)[u, u^{-1}]$.

%For a fixed vector bundle $\mc E$, we set
%\[
%S_{\mc E} = \lbrace \wedge_{-1} \left( 
%x_{1}^{m_{1}}
%\cdots x_{s}^{m_{s}} \mc E \right) 
%\mid (m_{1}, \ldots, m_{s})
%\in \Z^{s} \setminus \lbrace 0 \rbrace
%\rbrace \subset K_{\mb T}(M).
%\]

%%%%%%%%%%%%%%%%%%%%%%%%%%%%%%%%%%%%%%%%%%%%%%%%%%%%%%%%%%%%%%%%%%%%%

\subsection{Flag bundle}
\label{subsec:flag}
Let $\mc E$ be a $\mb T$-equivariant vector bundle 
over smooth $\mb T$-variety $B$. 
and consider a sequence of integers
$r \ge d_{1} \ge d_{2} \ge \cdots \ge d_{m} \ge 0$.
We consider the flag bundle of type $A_{m}$
\[
Fl(\mc E; d_{1}, d_{2}, \ldots, d_{m}) =
\left \lbrace \left.
( \mc E|_{b} \supset F_{1} \supset \cdots \supset F_{m} ) 
\ \right| \ b \in B, F_{k} \cong \C^{d_{k}} \ 
(k=1, \ldots, m)
\right \rbrace.
\]

We fix a non-negative integer $d$, and set 
$F=Fl(\mc E; d)$.
We write $\Pi^{\ast} \mc E$ by the same letter $\mc E$, 
and by $\mc V$ the universal sub-bundle on $F$, 
that is, 
the sub-bundle of $\mc E$ whose fibre
$\mc V_{\tilde{b}}$ over the point in $\tilde{b}
%=(V_{\tilde{b} \to \mc E|_{b}) 
\in F$ is equal to the subspace 
corresponding to $\tilde{b}$.
For a polynomial $g(Y_{1}, \ldots, Y_{d})
\in \kloc_{\mb T}(\mc M)[Y_{1}, \ldots, Y_{d}]$, 
we consider 
$K$-theory class $\psi=g(\wedge^{\bullet} \mc V)$ in 
$\kloc_{\mb T}
(\mc M)$.
Using a computational method inspired by \cite{M} and 
\cite{NY}, 
we compute the push-forward of 
$g ( \wedge^{\bullet} \mc V)$  
by the projection
\[
\Pi \colon
F=Fl(\mc E; d)
\to
B.
\]

Let us introduce the flag bundle $\Pi^{'} \colon
F^{'}=Fl(\mc E; d-1)$ of $(d-1)$-dimensional sub-spaces, 
and write by 
$\mc V^{'}$ the universal sub-bundle on $F'$ 
tentatively.
We take another torus $\C^{\ast}_{u}=\C^{\ast}$ with
the coordinate $u \in \C^{\ast}_{u}$ as above.
We consider the $\mb T \times \C^{\ast}_{u}$-equivariant 
line bundle $L = \mo_{F^{'}} \otimes \C^{\ast}_{u}$ over $F^{'}$, 
and the projective bundle
\[
\mc M
=
\PP_{F^{'}} \left( \left( \mc E / \mc V^{'}
\right) \oplus L \right)
\]
over $F^{'}$.
We call $\mc M$ the {\it master space}.
We consider a fibre-wise $\C^{\ast}_{u}$-action on $\mc M$
by $\C^{\ast}_{u}$-action on $L$ and the trivial action
on $\mc E / \mc V^{'}$.

For the fixed points set $\mc M^{\C^{\ast}_{u}}$,
we have a obvious decomposition
\[
\mc M^{\C^{\ast}_{u}}
=
\mc M_{+} \sqcup \mc M_{\exc}
\]  
of connected components, where
$\mc M_{+}=\PP_{F^{'}}\left( \mc E / \mc V^{'}
\right) = Fl(\mc E; d, d-1)$ and 
$\mc M_{\exc}=\PP_{F^{'}}(L)=F^{'}$.
Then we recover $F$ by the projection
$\mc M_{+} = Fl(\mc E; d, d-1) \to F=Fl(\mc E; d)$ 
forgetting the $d-1$-dimensional flags.
We write by the same letters $\mc V^{'}$ and $\mc V$ 
the pull-backs of $\mc V^{'}$ on $F^{'}$ to $\mc M$ and 
$\mc M_{+}$ and
$\mc V$ on $F$ to $\mc M_{+}$ respectively.

For the tautological bundle $\mo_{\mc M}(-1)$ of the
projective bundle $\mc M=\PP_{F^{'}} \left( \left(
 \mc E / \mc V^{'}
\right) \oplus L \right)$, we consider 
the inclusion of $\mo_{\mc M}(-1)$ 
into $\left( \mc E / \mc V^{'} \right) \oplus L$. 
We take the inverse image of 
$\mo_{\mc M}(-1)$ by the quotient map 
$\mc E \oplus L \to \mc E/\mc V^{'} \oplus L$, and 
write it by the same letter 
$\mc V$ as ones on $\mc M_{+}$ and $F$.
We have isomorphisms $\mc V|_{\mc M_{+}} \cong \mc V$ 
and $\mc V|_{\mc M_{\text{exc}}} \cong \mc V' \oplus L$, 
where $\mc V$ and $\mc V'$ have the trivial 
$\C^{\ast}_{u}$-actions.
We see that $\mo_{\mc M}(-1) = \mc V/\mc V'$ restricts 
to $\mo_{\mc M}(-1)|_{\mc M_{+}} = \mc V/\mc V^{'}$ on
$\mc M_{+}$ and 
$\mo_{\mc M}(-1)|_{\mc M_{\exc}} = L$ on 
$\mc M_{\exc}$.
Hence the normal bundles $N_{+}
\in K_{\mb T \times \C^{\ast}_{u}}(\mc M_{+})$ of
$\mc M_{+}$ and  
$N_{\exc}\in K_{\mb T \times \C^{\ast}_{u}}
(\mc M_{\exc})$ of 
$\mc M_{-}$ in $\mc M$ are 
written as
\begin{align}
\label{N+}
N_{+}
&=\mo_{\mc M}(1) \otimes L|_{\mc M_{+}} 
=
\C_{u} \otimes \mc V^{\vee} -
\C_{u} \otimes \mc V^{' \vee},
\\
\label{N-}
N_{\exc}
&=\mo_{\mc M}(1)|_{\mc M_{\exc}} \otimes 
\mc E / \mc V^{'}
=
\C_{u^{-1}} \otimes \mc E - 
\C_{u^{-1}} \otimes \mc V^{'}.
\end{align}
Furthermore we set 
$\Theta=\mc Hom(\mc V^{'}, \mo_{\mc M}(-1))$ 
on $\mc M$ so that we have
the restrictions 
\begin{align}
\label{Theta+}
\Theta|_{\mc M_{+}} 
&=\mc Hom(\mc V^{'}, \mc V/\mc V^{'}),\\
\label{Thetaexc}
\Theta|_{\mc M_{\exc}} 
&=\mc Hom(\mc V^{'}, \mo_{\mc M_{\exc}} \otimes \C_{u})
= \mc V^{' \vee} \otimes \C_{u}.
\end{align}
We write $\Theta|_{\mc M_{+}}$ in \eqref{Theta+} 
by $\Theta$.
This is equal to
the relative tangent bundle of the projection 
$\mc M
\to F$. 

Let $\iota_{+}$ and $\iota_{\exc}$ be
embeddings of 
$\mc M_{+}$ and $\mc M_{\exc}$ into $\mc M$.
We write by $\widehat{\Pi} \colon \mc M \to B$
the projection, and set
$\widetilde{\Pi}= \widehat{\Pi}|_{\mc M_{+}}$ and 
$\Pi^{'}= \widehat{\Pi}|_{\mc M_{\exc}}$.
By the similar argument as in the proof of 
\cite[Lemma 3.1]{N}, we see that 
$\wedge_{-1} (N_{+})^{\vee}$ and $\wedge_{-1} 
(N_{\exc})^{\vee}$ 
are invertible elements in 
\[
\kloc_{\mb T \times \C^{\ast}_{u}}(\mc M^{\C^{\ast}_{u}})
=\kloc_{\mb T \times \C^{\ast}_{u}}(\mc M_{+})
\oplus 
\kloc_{\mb T \times \C^{\ast}_{u}}(\mc M_{\exc}).
\]

We use the localization theorem
\cite[Theorem 3.2]{N} to prove the following lemma.
\begin{lem} 
\label{lem:loc}
We have the following commutative diagram
\[
\xymatrix{
\kloc_{\mb{T} \times \C^{\ast}_{u}} (\mc M) 
\ar[d]_{\widehat{\Pi}_{\ast} } \ar[r]^{\cong} & 
\kloc_{\mb{T} \times \C^{\ast}_{u}} 
(\mc M^{\C^{\ast}_{u}})  
\ar[d]^{\widetilde{\Pi}_{\ast} + \Pi^{'}_{\ast}} \\
\kloc_{\mb{T} \times \C^{\ast}_{u}} (B) 
\ar@{=}[r] & 
S_{\mb T \times \C^{\ast}_{u}, B}^{-1} K_{\mb{T}} (B)[u^{\pm 1}] 
}
\]
where the upper horizontal arrow is given by 
${\iota^{\ast}_{+} \over \wedge_{-1} N^{\vee}_{+}
 }+{\iota^{\ast}_{\exc} \over \wedge_{-1} 
 N^{\vee}_{\exc}
 }$.
\end{lem}
\proof
The argument is similar with slight modification 
because we treat
$\mc M^{\C^{\ast}_{u}}$ but not 
$\mc M^{\mb T \times \C^{\ast}_{u}}$.
\endproof

In particular, for any class $\psi \in K_{\mb T \times \C^{\ast}_{u}}
(\mc M)$, we have
\begin{align}
\label{localization}
\widehat{\Pi}_{\ast} \psi = \widetilde{\Pi}_{\ast} 
{\psi|_{\mc M_{+}} \over \wedge_{-1}(N_{+})^{\vee}} + 
\Pi^{'}_{\ast} {\psi|_{\mc M_{\exc}}
\over 
\wedge_{-1}(N_{\exc})^{\vee}} \in 
S_{\mb T \times \C^{\ast}_{u}, \mc M}^{-1} K_{\mb{T}} (B)[u^{\pm 1}]
\end{align}
by  Lemma \ref{lem:loc}.
%%%%%%%%%%%%%%%%%%%%%%%%%%%%%%%%%%%%%%%%%%%%%%%%%%%
%%%%%%%%%%%%%%%%%%%%%%%%%%%%%%%%%%%%%%%%%%%%%%%%%%%
%%%%%%%%%%%%%%%%%%%%%%%%%%%%%%%%%%%%%%%%%%%%%%%%%%%
\section{Wall-crossing formula}

To get the wall-crossing formula,
we use the equation \eqref{localization} in 
$\kloc_{\mb T}(B)$.
We introduce a map 
$\kloc_{\mb T}(\mc M^{\C^{\ast}_{u}}) 
=S_{\mb T \times \C^{\ast}_{u}, 
\mc M^{\C^{\ast}_{u}}}^{-1} 
K_{\mb{T}} (\mc M^{\C^{\ast}_{u}})[u^{\pm 1}] 
\to K_{\mb T}(\mc M^{\C^{\ast}_{u}})$ and
$\kloc_{\mb T}(B)[u^{\pm}] \to K_{\mb T}(B)$
sending $\varphi$ 
to a sum 
$\displaystyle \Res_{u = 0,\infty}
\varphi du$ of the residues at $u=0$ and $\infty$. 
They are compatible with the push-forwards 
$\widetilde{\Pi}_{\ast}$ and $\Pi^{'}_{\ast}$.
Since $\widehat{\Pi}_{\ast} \psi$ on the left hand side 
of \eqref{localization}
belongs to $K_{\mb T}(B)[u^{\pm 1}]$, we have
$\displaystyle \Res_{u = 0,\infty} \left( 
\widehat{\Pi}_{\ast} \psi \right) du=0$.
Hence we get
\begin{align}
\label{integral}
\Res_{u=0, \infty} \widetilde{\Pi}_{\ast} 
{\psi|_{\mc M_{+}} \over \wedge_{-1}(N_{+})^{\vee}} du=
-
\Res_{u=0, \infty} \Pi^{'}_{\ast} 
{\psi|_{\mc M_{\exc}}
\over 
\wedge_{-1}(N_{\exc})^{\vee}} du.
\end{align}

%%%%%%%%%%%%%%%%%%%%%%%%%%%%%%%%%%%%%%%%%%%%%%%%%%%%
\subsection{Residue formula}
For a polynomial $g(Y_{1}, \ldots, Y_{d}) \in 
\kloc_{\mb T}(\mc M)[Y_{1}, \ldots, Y_{d}]$,
we write the pull-backs of $g(Y_{1}, \ldots, Y_{d})$ by
the same letter. 

In $\kloc_{\mb T}(\mc M)$, we modify the $K$-theory class 
$\psi=g(\wedge^{\bullet} \mc V)$ as 
\[
\tilde{\psi} = {g(\wedge^{\bullet} \mc V) \cdot \wedge_{-1} \Theta^{\vee} 
\over (\mc V / \mc V^{'}) \cdot d } .
\]
By \eqref{N+}, we have 
\[
\Res_{u=0, \infty}
{\tilde{\psi}|_{\mc M_{+}} \over 
\wedge_{-1}(N_{+})^{\vee} } du
=
\Res_{u=0, \infty}
{\tilde{\psi}|_{\mc M_{+}} \over 1 - u^{-1} (\mc V 
/ \mc V^{'} )}du
=
- {g(\wedge^{\bullet} \mc V) \cdot \wedge_{-1} 
\Theta^{\vee} \over d}.
\]
Since $\Theta$ is the relative tangent bundle
of $\mc M_{+} = Fl(\E;d, d-1) \to F=Fl(\E;d)$, we have 
\begin{align}
\label{f(v)}
\displaystyle \Res_{u=0, \infty} 
\widetilde{\Pi}_{\ast} 
{\tilde{\psi}|_{\mc M_{+}} \over 
\wedge_{-1}(N_{+})^{\vee} } du 
=
- \Pi_{\ast} g(\wedge^{\bullet} \mc V).
\end{align}

On the other hand, by \eqref{N-} and \eqref{Thetaexc},
we can compute the right hand side
$- \displaystyle \Res_{u=0, \infty} \Pi^{'}_{\ast} 
{\tilde{\psi}|_{\mc M_{\exc}}
\over 
\wedge_{-1}(N_{\exc})^{\vee}} du$ in 
\eqref{integral}. 
As a result, we have 
\begin{align}
\label{euler}
\Pi_{\ast} g(\wedge^{\bullet} \mc V)
=
\Res_{u=0, \infty} \Pi^{'}_{\ast} 
\left[
g \left( \wedge^{\bullet} \left(
\mc V^{'} \oplus u \right) \right) \cdot 
{\wedge_{-1} \left( u \otimes \mc V^{'\vee} +
u^{-1} \otimes \mc V^{'} \right) 
\over
d \cdot \wedge_{-1} \left( u \otimes \mc E^{\vee}
\right)} \right]
{du \over u}.
\end{align}

%%%%%%%%%%%%%%%%%%%%%%%%%%%%%%%%%%%%%%%%%%%%%%%%%%%%
\subsection{Iterated residue}
We choose a polynomial $g(Y_{1},\ldots, Y_{d})
\in K_{\mb T}(\pt)[Y_{1}, \ldots, Y_{d}]$, and 
apply \eqref{euler} repeatedly to compute 
$\Pi_{\ast} g(\wedge^{\bullet} \mc V)$.

We introduce 
equivariant variables
$u_{1}, \ldots, u_{\ell}$, and set
\begin{align}
\nonumber
\psi_{\ell}(\mc E, \mc V)=
&g \left( \wedge^{\bullet} \left( 
\mc V \oplus \bigoplus_{i=1}^{\ell} u_{i} 
\right) \right) 
\cdot \prod_{1\le i \neq j \le \ell} (1 - u_{i}/u_{j})
\\
\label{itcoho}
&\cdot
\prod_{i=1}^{\ell} {\displaystyle
\wedge_{-1} (u_{i} \otimes \mc V^{\vee} + u_{i}^{-1} 
\otimes \mc V) 
\over
(d-i+1) \cdot \wedge_{-1} \left( u_{i} \otimes 
\mc E^{\vee} \right)
}
\in \kloc_{\mb T \times \C^{\ast}_{u_{1} } \times
\cdots \times \C^{\ast}_{u_{\ell}}} (Fl(\mc E; d- \ell)).
\end{align}
Note that we have a polynomial 
$g_{\ell}(Y_{1}, \ldots, Y_{d-\ell})$ such that 
$g_{\ell}(\wedge^{\bullet} \mc V) = 
\psi_{\ell} \left( \mc E, \mc V \right)$ on 
$Fl(\mc E; d- \ell)$.
In the next step, by \eqref{euler} applied to 
$g_{\ell}(\wedge^{\bullet} \mc V)$ in
$\kloc_{\mb T \times \C^{\ast}_{u_{1}} \times \cdots
\times \C^{\ast}_{u_{\ell}}} (\mc M)$
for $F=Fl(\mc E; d-\ell)$ and $F'=Fl(\mc E; d-\ell-1)$
replacing $\mb T$ with $\mb T \times \C^{\ast}_{u_{1}} 
\times \cdots \times \C^{\ast}_{u_{\ell}}$ and 
$u$ with $u_{\ell+1}$, 
we see that $\Pi_{\ast} \psi_{\ell}\left(
\mc E, \mc V \right)$ is equal to 
\begin{align*}
\Pi_{\ast} g_{\ell} \left( \wedge^{\bullet}
\mc V \right) 
&=
\Res_{u_{\ell+1}=0, \infty} \Pi_{\ast} 
\left[
g_{\ell} \left( \wedge^{\bullet} \left(
\mc V \oplus u_{\ell+1} \right) \right) \cdot 
{\wedge_{-1} \left( u_{\ell+1} \otimes \mc V^{\vee} 
\oplus u_{\ell+1}^{-1} \otimes \mc V \right) 
\over
d \cdot \wedge_{-1} \left( u_{\ell+1} \otimes 
\mc E^{\vee}
\right)} \right]
{du_{\ell+1} \over u_{\ell+1}}
\\
&=
\Res_{u_{\ell+1}=0, \infty} 
\Pi_{\ast} \psi_{\ell+1}(\mc E, \mc V)
{du_{\ell+1} \over u_{\ell+1}}. 
\end{align*}
We take $\ell$ times iterated residues
 $\displaystyle \Res_{u_{1}, \ldots, u_{\ell} =0, 
\infty}$, which means 
$\displaystyle \Res_{u_{1}=0, \infty} \circ \cdots \circ 
\Res_{u_{\ell}=0, \infty}$.
Then by induction on $\ell$, we have
\[
\Pi_{\ast} g(\wedge^{\bullet} \mc V) = 
\Res_{u_{1}, \ldots, u_{\ell} =0, 
\infty}
\Pi_{\ast} \psi_{\ell} (\mc E, \mc V)  
{du_{1} \cdots du_{\ell} 
\over 
u_{1} \cdots u_{\ell}}
\]
for any $\ell=0,\ldots, d$.
In particular, for $\ell=d$ we get the 
following theorem \cite[(3)]{WZ}.
\begin{thm}
\label{thm:main}
For the projection $\Pi \colon
F=Fl(\mc E; d)
\to
B$ of a flag bundle, we have
\begin{align}
\label{kresidue}
\Pi_{\ast} g(\wedge^{\bullet} \mc V)  
&=
\Res_{u_{1}, \ldots, u_{d} = 0, \infty}
\psi_{d} (\mc E) \cdot
{du_{1} \cdots du_{d} 
\over u_{1} \cdots u_{d}}
\in K_{\mb T} \left( B \right)
\end{align}
where $\psi_{d}(\mc E)=\psi_{d}(\mc E, 0)$ is a 
$K$-theory class 
defined in \eqref{itcoho}.
%in 
%$\kloc_{\mb T \times \C^{\ast}_{u_{1}} \times \cdots
%\times \C^{\ast}_{u_{d}} }(B)$
%defined by
Explicitly we can write
\begin{align}
\label{itresidue}
\psi_{d}(\mc E)=
{g \left( \wedge^{\bullet} ( u_{1} \oplus
\cdots \oplus u_{d} ) \right) \cdot
\prod_{1 \le i \neq j \le d} (1 - u_{i} / u_{j}) 
\over
d! \cdot \prod_{i=1}^{d} \wedge_{-1} \left( u_{i} \otimes
\mc E^{\vee} \right)}.
\end{align}
\end{thm}

By the similar argument using 
localization formula \cite[(1)]{GP} 
for the $\mb T$-equivariant Chow ring $CH^{\bullet}
_{\mb T} (M)$, 
we get the residue formula 
\cite[Corollary 1]{Z} for cohomology classes.
We re-write weight spaces by
$u_{1}=e^{z_{1}}, \ldots, u_{d}=e^{z_{d}}$, and
set $z_{1}=c_{1}(e^{\hbar_{1}}), \ldots, z_{d}
=c_{1}(e^{z_{d}})$ 
in $CH^{\bullet}_{\mb T}(\pt)$.
\begin{thm}
\label{thm:cohom}
For a cohomology class 
$g(c_{\bullet}(\mc V)) =g(c_{1}(\mathcal{V}), 
\ldots, c_{d}(\mathcal{V}))$ and 
the fundamental cycle $[F]$ of $F=Fl(\E; d)$,
we have
\begin{align}
\label{cohom}
  \Pi_{\ast} \left( 
  g(c_{\bullet} (\mc V)) \cap [F]
  \right) 
=  \Res_{z_{1} \cdots z_{d} = \infty}
{ f( z_{1}, \ldots, z_{d}) 
 \prod_{i \neq j} (z_{i} - z_{j}) 
 \over
 d! 
 \prod_{i=1}^{d} \Eu( e^{-z_{i}} \otimes 
 \mc E )}
d z_{1} \cdots d z_{d}
\in CH_{\bullet}^{\mb T}(B), 
\end{align}
where 
$f(x_{1}, \ldots, x_{d})=
g(e_{1}(x_{1}, \ldots, x_{d}), \ldots, 
e_{d}(x_{1}, \ldots, x_{d}))$, and
$\displaystyle 
\Res_{ z_{1} \cdots z_{d}=\infty}$ is 
taking iterations of residues 
$\displaystyle \Res_{ z_{1}=\infty} \circ \cdots 
\circ 
\Res_{ z_{d}=\infty}$. 
\end{thm}
If we write by $a_{1}, \ldots, a_{r}$ the Chern roots of
$\mc E$, then we hava
\begin{align}
\label{itcohom}
\Res_{z=\infty} {z^{k} dz \over 
\prod_{\alpha=1}^{r} - z +a_{\alpha}}
&=
\begin{cases}
(-1)^{r-1} h_{k+r-1}(a_{1}, \ldots, a_{r}) & k \ge r\\
0 & k < r.
\end{cases}
\end{align}
Here $h_{k}$ is the complete 
homogeneous symmetric function of degree $k$ 
\[
h_{k}(x_{1}, \ldots, x_{r} )=
\begin{cases}
\sum_{ \ell_{1}+ \cdots + \ell_{r} = k } 
x_{1}^{\ell_{1}} \cdots x_{r}^{\ell_{r}}& k \ge 0\\
0 & k < 0.
\end{cases}
\]
It can be checked directly as in \eqref{itclass}
below, and also by applying \eqref{cohom} for
$\PP(W)=Fl_{\pt}(W; 1)$ and $\psi(\mc V)=
c_{1}(\mc V)^{k}$.
Using \eqref{itcohom}, we can compute the linear map $\displaystyle 
\Res_{z_{1} \cdots z_{d}=\infty}$ in the right
hand side of \eqref{cohom}.

%%%%%%%%%%%%%%%%%%%%%%%%%%%%%%%%%%%%%%%%%%%%%%%%%%%%
%%%%%%%%%%%%%%%%%%%%%%%%%%%%%%%%%%%%%%%%%%%%%%%%%%%
%%%%%%%%%%%%%%%%%%%%%%%%%%%%%%%%%%%%%%%%%%%%%%%%%%%
%%%%%%%%%%%%%%%%%%%%%%%%%%%%%%%%%%%%%%%%%%%%%%%%%%%
\section{Application}
We apply Theorem \ref{thm:main} to the Grassmannian 
$G(W;d)$ of $d$-dimensional
sub-spaces in $W=\C^{r}$.
This is the flag manifold $Fl_{\pt}(W; d)$ 
for $X=\pt$, and $\mc E=W$
is regarded as a vector bundle of rank $r$ over the point 
$\pt$.
We set $\mc W=\mo_{G(W;d)} \otimes W$, and write by 
$\mc V$ the tautological sub-bundle of $\mc W$. 
We set $\int_{G(W;d)} \varphi=\Pi_{\ast} \varphi$
for a $K$-theory class $\varphi$ on $G(W;d)$.

For a fixed basis of $W=\C^{r}$, 
we consider the diagonal torus 
$\mb T=(\C^{\ast})^{r}$ of 
$\GL(W)$. 
We define $\mb T$-action on $G(W,;d)$ and 
$\mc M=\PP_{G(W;d-1)}( \mc W/ \mc V^{'} \oplus
\C_{u})$ via the natural $\GL(W)$-action on 
$G(W;d)$ and $\mc M$.
For a coordinate $(x_{1}, \ldots, x_{r})$ of 
$\mb T$, we consider weight spaces 
$\C_{x_{1}}, \ldots, \C_{x_{r}}$ 
with the weights $x_{1}, \ldots, x_{r}$ as in the 
previous section.

To compute $f_{d} (\mc W)$ in \eqref{itresidue}, 
we recall the computation of residues.
Using coordinate change $w=u^{-1}$ near $u=\infty$, 
we have
\begin{align}
\Res_{u=0, \infty} \frac{u^{k} } 
{\prod_{ \alpha=1}^{r} (1-  u/ x_{\alpha})} 
\frac{du}{u} 
&=
\Res_{u=0} u^{k-1}  
\prod_{ \alpha=1}^{r} (1+ u/ x_{\alpha} +u^{2} / 
x_{\alpha}^{2} + \cdots +) du
\notag
\\
&+(-1)^{r-1} \Res_{w=0} w^{-k+r} x_{1} \cdots x_{r} 
\prod_{ \alpha=1}^{r} (1+  w x_{\alpha} + w^{2} 
x_{\alpha}^{2} + \cdots) \frac{dw}{w} 
\notag
\\
&=
h_{-k} (1/x_{1}, \ldots, 1/x_{r}) + (-1)^{r-1} x_{1} 
\cdots x_{r} \cdot h_{k-r} (x_{1}, \ldots, x_{r})
\notag
\\
&=
\begin{cases}
h_{- k} (x_{1}^{-1}, \ldots, x_{r}^{-1})
& k \le 0\\
0 &  0 < k < r\\
(-1)^{r-1} \cdot x_{1} \cdots x_{r} \cdot 
h_{ k - r} (x_{1}, \ldots, x_{r})
& k \ge r.
\end{cases}
\label{itclass}
\end{align}

%%%%%%%%%%%%%%%%%%%%%%%%%%%%%%%%%%%%%%%%%%%%%%%%%%%
\subsection{Determinantal formula}

We consider a family $\lbrace g_{i}(u) 
\rbrace_{i \in \Z_{\ge 0}}$ 
of one variable Laurent polynomials.
We set
\[
f(u_{1}, \ldots, u_{r} ) = 
\frac{\det ( g_{i }( u_{j}))_{1\le i, j \le r}}
{ \prod_{1\le i < j \le r} (1/u_{i} - 1/u_{j})}.
\]
We prepare the following lemma.
(cf. in \cite[Theorem 1.1]{NNSY}) 
\begin{thm}
\label{det}
We have
\begin{align*}
{(-1)^{r(r-1)/2} f(x_{1}, \ldots, x_{r}) \over 
(x_{1} \cdots x_{r})^{r}}
=
\det \left( \Phi (g_{i}(u) u ^{j-1}) \right)
_{1 \le i, j \le r},
\end{align*}
where $\Phi \colon \C [u, u^{-1}] \to 
\C [x_{1}^{\pm 1}, \ldots, x_{r}^{\pm 1}]^{S_{r}}$ is a 
$\C$-linear map and $\Phi(u^{k})$ is defined by 
\begin{align}
\label{phi}
\Phi(u^{k})=
\begin{cases}
h_{- k } (x_{1}^{-1}, \ldots, x_{r}^{-1}) / 
(x_{1} \cdots x_{r}) & k \le 0\\
(-1)^{r-1}  \cdot h_{k - r} (x_{1}, \ldots, x_{r}) 
& k > 0. \\
\end{cases}
\end{align}
\end{thm}
\proof
We give another proof here for the convenience of 
the reader. 
When $d=r$, we have $G(W;r)=\Hom^{\text{inj}}(V, W)
/\GL(V)=\pt$, where 
$\Hom^{\text{inj}}(V, W)$ is the set of injective 
linear maps $V=\C^{r} \to W=\C^{r}$, that is, 
isomorphisms.
Then the tautological bundle $\mc V=
\Hom^{\text{inj}}(V, W) \times V / \GL(V)$ over 
the point
is isomorphic to $W$ as $\mb T$-representation via 
gauge transformations by $\GL(V)$.
Hence we can write 
$f(x_{1}, \ldots, x_{r} )=
\int_{G(W;r)} f(\mc V)$.
In Theorem \ref{thm:main}, we compute 
\[
f \left( u_{1}, \ldots, u_{r} \right) \cdot
\prod_{1 \le i \neq j \le r} (1 - u_{i} / u_{j})
=
\det \left( g_{i }( u_{j}) \right)_{1\le i, j \le r} \prod_{1\le i < j \le r} (u_{i} - u_{j})
\]
expanding $\det \left( g_{i }( u_{j}) 
\right)_{1\le i, j \le r}$.
We see that permuatations of variables 
$u_{1}, \ldots, u_{r}$ in \eqref{itresidue} does not 
give any change by \eqref{itclass}.
Hence we have
\begin{align*}
f(x_{1}, \ldots, x_{r} ) &= 
\Res_{u_{1}, \ldots, u_{r} = 0, \infty}
\frac{
\prod_{i=1}^{r} g_{i} (u_{i})   
\prod_{i<j} (u_{i} -u_{j})
}
{ \prod_{i=1}^{r} \prod_{\alpha=1}^{r} \left(1 -u_{i} / x_{\alpha} \right)}
\cdot
\frac{du_{1} \cdots du_{r}}{u_{1} \cdots u_{r}}.
\end{align*}
We further expand $\prod_{i<j} (u_{i} -u_{j})$ to get
\begin{align*}
f(x_{1}, \ldots, x_{r} )
&= 
(-1)^{r(r-1)/2} 
 \sum_{\sigma \in \mk S_{r}}
\text{sign}(\sigma)
\Res_{u_{1}, \ldots, u_{r} = 0, \infty}
\prod_{i=1}^{r} 
\frac{ g_{i} (u_{i})
u_{i}^{\sigma(i) - 1}}
{\prod_{\alpha=1}^{r} \left(1 -u_{i} / x_{\alpha} \right)}
\cdot
\frac{du_{1} \cdots du_{r}}{u_{1} \cdots u_{r}}
\notag\\
&= 
(-1)^{r(r-1)/2} 
\det \left(
\Res_{u = 0, \infty}
\frac{ 
g_{i} (u)
u^{j - 1}}
{ \prod_{\alpha=1}^{r} \left(1 -u / x_{\alpha} \right)}
\cdot
\frac{du}{u }
\right).
\end{align*}
By \eqref{itclass}, we get the assertion.
\endproof

In particular, for a partition $\lambda=(\lambda_{1}, 
\lambda_{2}, \ldots, )$ with the length $\ell \le r$, 
we set
$g_{i}(u) = u^{\lambda_{i} + r- i+1 }$.  
Then $f(x_{1}, \ldots, x_{r})$ is equal to 
$(x_{1} \cdots x_{r})^{r} S_{\lambda}$ where
$S_{\lambda}(x_{1}, \ldots, x_{r})$ is the Schur polynomial.
In this case, Theorem \ref{det} gives the Jacobi-Trudi 
formula.

%We have 
%\begin{align*}
%S_{\lambda}(x_{1}, \ldots, x_{r}) 
%&= 
%\det_{1\le i, j \le r} ( 
%h_{\lambda_{i} - i + j}(x_{1}, \ldots, x_{r})).
%\end{align*}

%%%%%%%%%%%%%%%%%%%%%%%%%%%%%%%%%%%%%%%%%%%%%%%%%%%%%%%%%%%%%%
\subsection{Refined shifted factorial Grothendieck polynomial}
We prepare an infinite number of parameters
$b=(b_{1}, b_{2}, \ldots, )$, $\alpha=(\alpha_{1}, \alpha_{2}, \ldots, )$ and $\beta=(\beta_{0}, \beta_{1} ,\beta_{2}, \ldots, )$.
We consider the definition of the factorial 
Grothendieck polynomial by \cite[(2.12), (2.13)]{IN} and generalize it including the
refined canonical Grothendieck polynomial introduced in \cite{HJKSS} as follows.

We take a partition $\lambda=(\lambda_{1}, \ldots, \lambda_{r})$ 
where $\lambda_{1} \ge \cdots \ge \lambda_{r} \ge 0$. 
For $i=1, \ldots, r$, we consider polynomials 
\[
g_{i}(u) = [u|b]^{\lambda_{i} +r- i} \frac{(1 - \beta_{1} u) \cdots (1 - \beta_{i-1} u)}{ (1 - \alpha_{1} u) \cdots (1 - \alpha_{\lambda_{i}} u)} u
\] 
where $[u|b]^{k}=(u+b_{1}+\beta_{0} u b_{1}) \cdots 
(u +b_{k}+\beta_{0} u b_{k})$. 

\begin{defn}
\label{rsfg}
We call the polynomial
\begin{align*}
 G_{\lambda} (u_{1}, \ldots, u_{r} | b, \alpha, \beta) 
&=
{\det \left( g_{i} ( u_{j}) / u_{j} \right)
_{1\le i, j \le r} 
\over 
\prod_{1 \le i<j\le r} ( u_{i} - u_{j}) }
\end{align*}
refined shifted factorial Grothendieck polynomial.
\end{defn}
When $b=(0, 0, \ldots, )$, $\alpha=(0,0, \ldots, ), \beta_{1} = \cdots = \beta_{r} = - \beta_{0}$, 
the determinantal formula
\begin{align}
G_{\lambda} (x_{1}, \ldots, x_{r} |0, 0, \beta)
&=
\det
\left(
\sum_{m=0}^{i-1} \begin{pmatrix} i-1 \\ m \end{pmatrix} \beta_{0}^{m} h_{\lambda_{j} - i +j +m}(x_{1}, \ldots, x_{r})
\right)
\label{det2}
\end{align}
is obtained in \cite{L}.
When we further set $\beta_{0}=-1$ and 
$x_{i} = 1 - 1 / \alpha_{i}$, then
Definition \ref{rsfg} coincides with the traditional 
combinatorial definition of
the Grothendieck polynomial by \cite[Proposition 3.3]{RS} 
where the double Grothendieck polynomials are studied
(see also the references therein).

In \cite[3.11]{HIMN}, the generalization of \eqref{det2} is obtained
for arbitrary $b=(b_{1}, b_{2}, \ldots, )$.
See also \cite[Remark 5.4 (2)]{NN1}.
The special case is also obtained by \cite[Proposition 3.8]{I}, and the skew version by \cite[Corollary 4.2]{IwaoMS2022} 
using free-fermions.
On the other hand, in \cite[Theorem 1.3]{HJKSS} the generalization of \eqref{det2} to the case where $b=(0,0,\ldots)$ and
arbitrary $\alpha, \beta$ is obtained.
We generalize these formula to the case for arbitrary 
$b, \alpha, \beta$, 
applying Theorem \ref{det} to 
$\frac{(-1)^{r(r-1)/2} f(x_{1}, \ldots, x_{r})}{ (e_{1} \cdots e_{r})^{r}}=
G_{\lambda} (x_{1}, \ldots, x_{r} |b, \alpha, \beta)$
as follows.
\begin{thm}
\label{thm:groth}
We have
\begin{align*}
G_{\lambda} (x_{1}, \ldots, x_{r} |b, \alpha, \beta)
=
\det \left(\Phi \left( 
[u|b]^{\lambda_{i} +r- i} \frac{(1 - \beta_{1} u) \cdots (1 - \beta_{i-1} u)}
{ (1 - \alpha_{1} u) \cdots (1 - \alpha_{\lambda_{i}} u)} 
u ^{j} \right)
\right)_{1 \le i, j \le r},
\end{align*}
where $\Phi$ is a $\C$-linear map defined by \eqref{phi}.
\end{thm}

When $\alpha=0$, 
we remark that in the determinantal formula 
\cite[3.11]{HIMN}, components of the matrix
are written by linear combinations of Schubert classes.
In our formula, we do not know such a geometric 
interpretation in Thorem \ref{thm:groth} at the present.
But the matrix element in Theorem \ref{thm:groth} is
 a finite sum.
\\

%%%%%%%%%%%%%%%%%%%%%%%%%%%%%%%%%%%%%%%%%%%%%%%%
%%%%%%%%%%%%%%%%%%%%%%%%%%%%%%%%%%%%%%%%%%%%%%%%
%%%%%%%%%%%%%%%%%%%%%%%%%%%%%%%%%%%%%%%%%%%%%%%%
%%%%%%%%%%%%%%%%%%%%%%%%%%%%%%%%%%%%%%%%%%%%%%%%
%%%%%%%%%%%%%%%%%%%%%%%%%%%%%%%%%%%%%%%%%%%%%%%%

\noindent {\bf Funding and/or Conflicts of 
interests/Competing interests}

He is partially supported by Grant-in-Aid for Scientific 
Research 21K03180 and 17H06127, JSPS.
This work was partly supported by Osaka Central Advanced Mathematical
Institute: MEXT Joint Usage/Research Center on Mathematics and
Theoretical Physics JPMXP0619217849, and by the Research Institute for Mathematical Sciences,
an International Joint Usage/Research Center located in Kyoto University.
The authors have no conflicts of interest directly relevant to the content of this article.
Data sharing not applicable to this article as no datasets were generated or analysed during the current study.

%%===========================================================================================%%
%% If you are submitting to one of the Nature Portfolio journals, using the eJP submission   %%
%% system, please include the references within the manuscript file itself. You may do this  %%
%% by copying the reference list from your .bbl file, paste it into the main manuscript .tex %%
%% file, and delete the associated \verb+\bibliography+ commands.                            %%
%%===========================================================================================%%

\bibliography{residue}% common bib file
%% if required, the content of .bbl file can be included here once bbl is generated
%%\input sn-article.bbl

\end{document}